\documentclass[12pt]{amsart}
\usepackage{graphics}
\usepackage{amssymb}
\usepackage{verbatim}
\usepackage{amsmath}
\usepackage{amscd}
\theoremstyle{plain}
\newtheorem{theorem}{Theorem}[section]
\newtheorem{corollary}[theorem]{Corollary}
\newtheorem{proposition}[theorem]{Proposition}
\newtheorem{lemma}[theorem]{Lemma}

\theoremstyle{definition}
\newtheorem{definition}{Definition}[section]
\newtheorem{example}{Example}[section]

\theoremstyle{remark}
\newtheorem{remark}{Remark}

\newcommand{\cstar}{\ensuremath{\text{C}^{*}}\nobreakdash-\hspace{0pt}}

\newcommand{\dirlim}{\varinjlim}
\newcommand{\h}{\hat}
\newcommand{\x}{\mathbf{x}}
\newcommand{\y}{\mathbf{y}}
\newcommand{\fA}{\mathfrak{A}}
\newcommand{\fB}{\mathfrak{B}}
\newcommand{\fD}{\mathfrak{D}}
\newcommand{\cA}{\mathcal{A}}
\newcommand{\cC}{\mathcal{C}}
\newcommand{\cD}{\mathcal{D}}
\newcommand{\cE}{\mathcal{E}}
\newcommand{\cF}{\mathcal{F}}

\newcommand{\cH}{\mathcal{H}}
\newcommand{\cN}{\mathcal{N}}
\newcommand{\cO}{\mathcal{O}}
\newcommand{\cP}{\mathcal{P}}
\newcommand{\cQ}{\mathcal{Q}}
\newcommand{\cR}{\mathcal{R}}

\newcommand{\cV}{\mathcal{V}}
\newcommand{\cW}{\mathcal{W}}
\newcommand{\bC}{\mathbb{C}}
\newcommand{\bN}{\mathbb{N}}
\newcommand{\bR}{\mathbb{R}}
\newcommand{\bT}{\mathbb{T}}
\newcommand{\bZ}{\mathbb{Z}}
\newcommand{\al}{\alpha}
\newcommand{\p}{\varphi}

\newcommand{\ti}{\tilde}
\newcommand{\om}{\omega}
\newcommand{\De}{\Delta}
\newcommand{\G}{\Gamma}
\newcommand{\ba}{\backslash}

\newcommand{\Bs}{Bratteli system}
\newcommand{\ra}{\rightarrow}
\newcommand{\Ra}{\Rightarrow}
\newcommand{\ve}{(\mathcal{V},\mathcal{E})}

\newcommand{\xmax}{X_{\max}}
\newcommand{\xmin}{X_{\min}}
\newcommand{\xmaxc}{X\ba \xmax}
\newcommand{\xminc}{X\ba \xmin}
\newcommand{\dom}{\mbox{dom }}
\newcommand{\pds}{partial dynamical system }

\begin{document}
\addtolength{\textwidth}{25 pt}
% topmatter
\title[]{$\bZ$-analytic TAF algebras and partial dynamical systems}

\author{Justin~R. Peters} \author{Yiu T. Poon} \thanks{The first
author was supported in part by NSF grant DMS9504359}
\address{Department
of Mathematics\\
	Iowa State University\\
	Ames, IA, USA}
\email{peters@iastate.edu, ytpoon@iastate.edu}

% \keywords{keywords}
% \subjclass{Primary: subject; Secondary: subject} % \date{\date}

% abstract
\begin{abstract}
Connections between partial dynamical systems, a generalized notion of
partial dynamical systems defined by nested sequences of partial
hoeomorphisms, and triangular AF- algebras which admit an
integer-valued cocycle are established.
\end{abstract}
\maketitle

\section{Introduction}\label{s:intro}

Among the analytic triangular AF (TAF) algebras, the $\bZ$-analytic
algebras are those admitting an integer-valued cocycle. The purpose
of this paper is to explore the connection between these algebras and
dynamical systems.

The work of this paper may be, in a sense, considered a parallel
version for nonselfadjoint algebras to the program of \cite{hps92},
\cite{gps95} studying the connections between zero-dimensional
dynamical systems, crossed products, and their K-theoretic
invariants. In some respects, the situation here is simpler:
isomorhism of standard $\bZ$-analytic dynamical systems is equivalent
to conjugacy of the dynamical systems, so notions such as strong
orbit equivalence are not needed. Also, since the K-theory for
nonseladjoint algebras is not well developed, we do not attempt to
discuss K-theoretic invariants. On the other hand, there are
complications in the nonselfadjoint situation not encountered in the
self-adjoint case. First of all, the dynamical systems are only
partially defined; worse yet, the `dynamical systems' may require an
infinite sequence of partial homeomorphisms whose graphs are nested.
The algebras to which the dynamical systems are associated are
nonselfadjoint subalgebras of groupoid C$^*$-algebras (cf
\cite{ms89}), or may also be viewed as nonselfadjoint subalgebras of
generalized crossed products in the sense of Exel (\cite{re94},
\cite{re95}).

By a \emph{spectral triple} $(X, \cP, \cR)$ we mean $X$ is a compact
space, which in our setting will always be zero-dimensional, $\cR$ is
a groupoid which is an equivalence relation on $X$ having unit space
$X,$ and $\cP$ an open subset of $\cR$ with $\cP \cup \cP^{-1} = \cR,
\quad \cP \circ \cP \subseteq \cP$, and $\cP \cap \cP^{-1} = X$
\cite{jr80}. $\cP$ determines a nonselfadjoint algebra $\fA = \fA(\cP)
\subset \text{C}^*(\cR)$ \cite{ms89}, which is `triangular' in the
sense that $\fA \cap \fA^* \cong C(X) \subset \text{C}^*(\cR).$ The
spectral triples we will study have the property that for all $(x, y)
\in \cP$ the set \[ \{ z: (x, z),\ (z, y) \in \cP \}\] has finite
cardinality.  In other words, regarding $\cP$ as giving an ordering to
the equivalence classes of $X,$ the order type is that of a subset of
the integers.

For example, if $X$ is a compact metrix space, $\p: X \to X$ a
homeomorphism, set \[ \cR = \{ (x, \p^n(x)): \ x\in X,\ n \in \bZ \},\]
\[ \cP = \{
(x, \p^n(x)): \ x\in X,\  n \geq 0 \},\] then C$^*(\cR)$ is
the crossed product algebra $C(X)\times_{\p}\bZ, \quad$
$\fA(\cP)$ the semicrossed product algebra $C(X)\times_{\p}\bZ^{+}$.
Suppose $X$ is zero-dimensional, and $\p$ a partial homeomorphism on
$X.$ That is, $\p$ is a homeomorphism from $Dom(\p)$ to $Ran(\p),$
where both $Dom(\p),\ Ran(\p)$ are  open  subsets of $X.$
Then, under certain conditions,  C$^*(\cR)$ is an AF algebra, and $\fA(\cP)$ a standard
$\bZ$-analytic subalgebra.  The partial homomorphism may arise as the
restriction of a homeomorphism of $X,$ (the case studied in \cite{hps92}) or it may not admit an
extension to a homeomorphism.

Returning to our general setting, suppose $(X, \cP, \cR)$ is a
spectral triple such that each equivalence classes of $X$ has the
order type of a subset of the integers. In that case it is possible
to define a partial mapping $\Phi$ on $X$ which maps each point to
its immediate successor in the ordering. $Dom(\Phi)$ is the set of
all points which have a successor. In general, such a map need not be
continuous. Spectral triples of this kind can arise from nested
sequences of partial homoemorphisms, $\{\p_n\}_{n\in \bZ}.$ By
\emph{nested} sequence we mean that $\p_{-n} = \p_n^{-1} (n\in \bZ),\
\p_0$ is the identity map, and that the graphs are nested in the
sense \[ \G(\p_n\circ \p_m)\subset \G(\p_{n+m}),\quad n, m \in \bZ .\]

Given a partial homeomorphism $\p$ of a zero dimensional space $X,$
under what conditions is the groupoid
$\cR(X, \p) := \{ (x, \p^n(x)) : x, \in X,\ n \in \bZ\}$ an AF
groupoid? We prove that this happens exactly when $(X, \p)$ is conjugate to $(B,
\psi),$ where $B$ is an ordered Bratteli compactum and $\psi$ the
Ver\v{s}ik map on $B.$ And this is the case if and only if for any
clopen set $U$ containing the complement of the range of $\p,$ each
$x \in X$ belongs to the forward orbit of some point $u \in U.$ (Cor.~\ref{c:afgpd})

The more general case in which $\cR$ is the union of the graphs
$\G(\p_n)$ of a nested sequence of partial homeomorphims is precisely
the case in which $\cR$ admits a (continuous) integer valued cocycle.
Here, too, we have necessary and sufficient conditions, given in
terms of dynamical systems, for $
\cR$ to be an AF groupoid.

Though our goal initially
was to study nonselfadjoint algebras, our approach leads to some new
results in AF algebras and dynamical systems.  Our approach, however,
comes naturally from looking at nonselfadjoint algebras. The nonselfadjoint
algebras we study are the strongly maximal subalgebras of AF algebras
which admit an integer-valued cocycle, i.e., the $\bZ$-analytic algebras.
Theorem \ref{t:zanal} in section \ref{s:char} gives various
characterizations of $\bZ$-analytic algebras (and spectral triples)
assuming the enveloping groupoid is AF, and an analogous result for
standard $\bZ$-analytic algebras is presented in \ref{t:stdzanal}. In
section \ref{s:pds} two standard $\bZ$-analytic algebras are shown to
be isomorphic iff their partial dynamical systems are conjugate.
(There is a corresponding result for semicrossed products and
dynamical systems; cf. \cite{scp92a}). We return to this theme again
in the context of nested sequences of partial homeomorphisms in
section \ref{s:nested} in which the dynamics may not be continuous
(Theorem \ref{t:nconj}). An example of a $\bZ$-analytic semigroupoid which is
not standard is given in Theorem \ref{e:dh}.  This example admits a
semi-saturation, and in the final subsection \ref{ss:sat} we consider
necessary and sufficient conditions for a nested sequence of partial
homeomorphisms to admit a semi-saturation.  Example \ref{e:nonsemisat}
is one that does not admit a semi-saturation.

\section{Preliminaries} \label{s:prelim}

We begin this section with a review of terminology for AF groupoids
and semi-groupoids.
\cite{scp92bk}.

Let $\fA$ be a triangular AF algebra, or simply TAF algebra, with
canonical masa $\fD = \fA \cap \fA^*,$ and \cstar envelope $\fB .$ In
the spectral triple for $\fA,$ denoted $ \ (X, \cP, \cR),$ $X$ is the
Gelfand space of the commutative \cstar algebra $\fD, \ \cR $ is the
AF groupoid of $\fB,$ and $\cP$ is the semigroupoid corresponding to
the subalgebra $\fA.$ If $v$ is a matrix unit of $\fA,$ or more
generally if $v$ is a $\fD$ - normalizing partial isometry in $\fA,$
let $\h{v}$ denote its support: i.e., $\h{v}$ is the support set of
$v,$ viewed as a function on the groupoid $\cR.$ Thus,
\begin{equation*}
\cP = \cup \{ \h{v} : v \text{ is a matrix unit of } \fA \} .
\end{equation*}
The sets $\h{v},$ as $v$ ranges over the matrix units of $\fA,$ form
a basis for the topology of $\cP .$ (And the sets $\h{v},$ form a
basis for the topology of $\cR,$ as $v$ ranges over the matrix units
of $\fB .$ )

A \emph{cocycle} (or, more precisely, a real-valued 1-cocycle) on
$\cR$ is a map $c: \cR \to \bR $ satisfying \begin{itemize} \item the
cocycle condition: for all $(x, y),\ (y, z) \in \cR, \ c(x, y) + c(y,
z) = c(x, z) $

\item continuity: $c$ is continuous from $\cR$ to $\bR.$

\end{itemize}

A TAF algebra $\fA$ with spectral triple $(X, \cP, \cR)$ is said to
be \emph{analytic} if there is a (real-valued 1-) cocycle $c$ such
that $c^{-1}([0, \infty)) = \cP $.
%\item $c^{-1}(0) = \Delta(X) := \{ (x, x): x \in X \} .$
We say in this case that $\fA$ is the analytic TAF algebra defined by
$c,$ and write $\fA = \fA_c .$ $\fA$ is called $\bZ$-analytic if the
cocycle $c$ can be chosen to be integer valued.

There is a
proper subclass of the $\bZ$-analytic algebras, called the
\emph{standard} $\bZ$-analytic algebras. If $\fA$ is any strongly
maximal TAF
subalgebra of an AF algebra $\fB,$ then by \cite{ppw90} Lemma 1.1
there
is a sequence $\{ \fB_n \}$ of finite dimensional \cstar algebras of
$\fB$, say $\fB_n = \bigoplus _{k=1}^{r(n)} M_{m(n, k)}$, and a set
of matrix units $\{ e_{ij}^{(nk)} \}$ for $\cup _{n=1}^{\infty}
\fB_n$ such that
\[ \fA_n:= \fA\cap \fB_n = \bigoplus _{k=1}^{r(n)} T_{m(n, k)}, \]
where $T_{m(n, k)}$ is the upper triangular subalgebra of $ M_{m(n,
k)},$ and $\fA = \dirlim \fA_n$ and $\fB = \dirlim \fB_n.$ Then $\fA$
is called
\emph{standard} if the embeddings $\fA_n \to \fA_{n+1}$ can be chosen
to be standard \cite{powa93};
i.e., standard embeddings in the sense of Effros (\cite{bb98}).
Since $\bZ$-analytic algebras are defined by means of the existence
of an integer-valued cocycle,
it is not obvious from the definition that standard $\bZ$-analytic
algebras form a subclass of the
$\bZ$-analytic. If $\fA$ is $\bZ$-analytic with spectral triple $(X,
\cP, \cR),$
and hence admits an integer-valued cocycle $d,$ it follows that if $x
< y$ are two points in $X$
belonging to the same equivalence class, then there are at most
finitely many points in the equivalence
class between $x$ and $y:$ indeed, the number of such points can be
at most $d(x, y).$ Thus one can define
an integer-valued functin on $\h{d},$ on $\cP,$ called the counting
cocycle:
\[ \h{d}(x, y) = 1 + \text{ the number of points in the orbit between
} x \text{ and } y \, .\]
It is clear that $\h{d}$ satisfies the additivity property for
cocycles; however in general it is not continuous. From \cite{powa93}
it is known that the standard $\bZ$-analytic algebras are precisely
the ones for which the counting cocycle is continuous. (cf
\ref{t:stdzanal})

\begin{definition} \label{d:pds} A partial dynamical system is a
quadruple $(X, \xmax, \xmin, \p)$ where $X$ is a compact metric
space, $\xmax,\ \xmin$ are closed subsets, and $\p: X\ba\xmax \to
X\ba\xmin$ is a homeomorphism. $\p$ is called a partial homeomorphism
on $X.$ The graph of $\p$ will be denoted $\G_p.$ \end{definition}

In this paper $X$ will always be assumed to be zero dimensional. For
$n$ a positive integer, $\p^n$ will denote the $n$-fold composite of
$\p,$ and for $n$ a negative integer, $\p^n$ will denote the $n$-fold
composite of $\p^{-1}.$ Also we adhere to the convention that $\p^0 =
id_X.$

For $U \subset X,$ we write $\p(U)$ to mean $\p(U\ba\xmax),$ and
similarly $\p^n(U) = \p^n(U\cap \text{dom}(\p^n)).$ Since $\xmax$ is
closed it follows that if $U$ is open, $\p(U)$ is open; and $\p^n(U)$
is open ($n \in \bZ$).

\section{Characterizations of $\bZ$-analytic and standard
$\bZ$-analytic algebras} \label{s:char}

We begin with several characterizations of $\bZ$-analytic and
standard $\bZ$-analytic algebras.
For subsets $A, B$ of a operator algebra $\fA,$ we will write
$[A\cdot B]$ to denote the closed linear span of the set $\{ab: \ a
\in A,\ b \in B \}.$ Let $\bT$ be the circle formed by identifying the
end points of $[0,\ 1]$.

\begin{theorem} \label{t:zanal} Let $\fA$ be a strongly maximal TAF
algebra with spectral triple $(X, \cP, \cR).$ The following conditions are
equivalent: \\ \begin{enumerate} \item \label{i1} C$^*(\fA)$ admits a
continuous action $\al$ of $\bT$, with fixed point algebra the
diagonal $\fD,$ and such that
$$\fA = \{ b \in \text{C}^*(\fA):\
\int_{\bT} \al_t(b) e^{2\pi int}\,dt = 0 \},$$ for $n \ge 1.$

\item \label{i2} There is a sequence $\{ \fA(n) \}, n \ge 0, $ of
closed linear subspaces of $\fA$ satisfying: \begin{enumerate} \item
$\fA(0)$ is the diagonal $\fD$ \label{i2a} \item $\fA(n) \cap \fA(m)
= (0), \ n \neq m$ \label{i2b} \item $[\fA(n)\cdot \fA(m)] \subset
\fA(n+m),$ for all $n, m \ge 0$ \label{i2c} \item
$\oplus_{n=0}^{\infty} \fA(n)$ is dense in $\fA$ \label{i2d}
\end{enumerate}

\item \label{i3} $ \fA$ is $\bZ$-analytic

\item \label{i4} There is a nested sequence $\{\p_n \}$ of partial
homeomorphisms of $X,$ with $\p_0 = id_X,$ the graph $\G_{\p_n}$ open
in $\cP,\ (n \ge 0)$ and such that \[ \cP = \sqcup_{n=0}^{\infty}
\G_{\p_n} \quad \text{(disjoint union)} \] \end{enumerate}
\end{theorem}

\begin{proof}
\ref{i1} $\Longrightarrow$ \ref{i2}.
Note that $\al$ acts on $\fA.$
Set $\fA(n) := \{ a \in \fA : a
= \int_0^1 \al_t(b)e^{-2\pi int}\,dt \ , \text{ for some } b \in \fA
\}$, for $n\ge 0$ and $\fA(n) :=(0)$, otherwise. First we show
\ref{i2a}. Since
the fixed point algebra of $\al$ is the diagonal $\fD,$ clearly
$\fA(0) \supset \fD.$ Now let $a \in \fA(0)$; then, if $ a = \int_0^1
\al_t(b)\,dt$
%\begin{multline*}
$$\al_s(a) = \int_0^1
\al_s(\al_t(b))\,dt = \int_0^1 \al_{s+t}(b)\,dt = \int_0^1
\al_t(b)\,dt = a,
%\end{multline*}
$$
where the last equality uses the translation invariance of Lebesgue
measure on $\bT.$ Thus, the fixed point algebra of $\al,$ which by
assumption is $\fD,$ contains $\fA(0).$ Thus $\fA(0) = \fD.$

Next, observe that $\fA(n) = \{ a \in \fA:\ \al_t(a) = e^{2\pi
int}a,\ t \in \bT\}.$ Indeed, let \mbox{$a = \int_0^1
\al_t(b)e^{-2\pi int}\, dt \ \ (b\in \fA).$} Then, for $s\in \bT,$
%\begin{multline*}
$$\al_s(a) = \int_0^1 \al_{s+t}(b)e^{-2\pi int}\, dt = \int_0^1
\al_u(b) e^{-2\pi in(u-s)}\, du = (e^{2\pi is})^n a. $$
%\end{multline*}
Conversely, if $\al_s(a) = (e^{2\pi is})^n a,$ a similar calculation
shows that $ a \in \fA(n).$ Thus, \ref{i2b} is clear, as is
\ref{i2c}.

To show \ref{i2d}, let $\om $ be a continuous linear functional on
$\fA,$ and suppose
\[ <a, \om> = 0 \ \text{ for } a \in \fA(n), \ n =0, 1, 2,\dots. \]
Fix $b \in \fA,$ and consider the continuous function $f$ on $\bT, \
f(t) = <\al_t(b), \om>.$ By the characterization of $\fA(n)$ above,
we have $\h{f}(n) = 0$ for all $n.$ Thus, $
f = 0,$ and hence $<b, \om> = 0.$ Since $b \in \fA$ was arbitrary, it
follows that $\om = 0.$ Thus, $\oplus_{n\in \bZ} \fA(n)$ is dense in
$\fA.$ Since by assumption the spaces $\fA(-n) \ (n>0)$ are $(0),$ we
have that $\oplus_{n=0}^{\infty}\fA(n)$ is dense in in $\fA.$

\ref{i2} $\Longrightarrow$ \ref{i3} Condition \ref{i2c} of \ref{i2}
implies that each $\fA(n)$ is a closed $\fD-$bimodule. Hence, by
\cite{ppw90} Theorem 2.2, $\fA(n)$ is spanned by the matrix units it
contains. Write
\[ \h{\fA}(n) = \cup\{ \h{v}: \ v \text{ a matrix unit in } \fA(n)\}.
\]
We claim: $\cup_{n=0}^{\infty} \h{\fA}(n) = \cP.$

Let $(x_0, y_0) \in \cP,$ and suppose $(x_0, y_0) \notin \h{\fA}(n),
\ n \geq 0.$ Define a representation $\pi$ of $\fA$ as follows: let
$\cH_{\pi}$ be a Hilbert space with orthonormal basis $\{\xi_x\}_{x
\in \cO(x_0)},$ and define $\pi$ on matrix units by \begin{equation*}
\pi(v)\xi_x =
\begin{cases}
\xi_y \text{ if (x, y) } \in \h{v}, \ x \in \cO(x_0) \\ 0 \text{
otherwise.}
\end{cases}
\end{equation*}
By \cite{jojp95} $\pi$ extends to a representation of $\fA.$ By
supposition, $(\pi(a)\xi_{x_0}, \xi_{y_0}) = 0 $ for $a \in \fA(n),$
hence for $a \in \oplus_{n=0}^{\infty} \fA(n).$ Since
$\oplus_{n=0}^{\infty} \fA(n)$ is dense in $\fA,$ we have that
$(\pi(a)\xi_{x_0}, \xi_{y_0}) = 0 $ for $a \in \fA.$ That is
impossible, since there is a matrix unit $v \in \fA$ with $(x_0, y_0)
\in \h{v},$ hence $(\pi(v)\xi_{x_0}, \xi_{y_0}) = 1.$ This proves the
claim.

Next, $\h{\fA}(n) \cap \h{\fA}(m) = \emptyset, $ for $n \neq m.$
Indeed, the intersection $\h{\fA}(n) \cap \h{\fA}(m)$ is open, so if
it is nonempty there is a matrix unit $v$ such that $\h{v}$ lies in
the intersection. But then $v \in \fA(n) \cap \fA(m),$ contradicting
\ref{i2b}.

Define a cocycle $d: \cP \to \bZ,$ by $d(x, y) = n$ if $(x, y) =
\h{\fA}(n).$ Then $d$ is well-defined, and continuous since
$d^{-1}(n) = \h{\fA}(n)$ is open in $\cP.$ Also, if $(x, y) \in
\h{\fA}(n),\ (y, z) \in \h{\fA}(m),$ then by \ref{i2c} $(x, z) \in
\h{\fA}(n+m),$ so $d(x, y) + d(y, z) = d(x, z),$ i.e., the cocycle
condition is satisfied on $\cP.$ $d$ can be extended to a cocycle on
$\cR$ by setting $d(y, x) = -n$ if $d(x, y) = n.$ One easily checks
that the cocycle condition is satisfied on the groupoid $\cR.$

\ref{i3} $\Longrightarrow$ \ref{i4} Let $\p_n$ be the partial
homeomorphism on $X$ whose graph $\G_{\p_n}$ is $d^{-1}(n).$ It is
clear from \ref{i2c} that  $\{\p_n \}$ form a nested sequence and
that the graphs $\G_{\p_n}$ have the required properties.

\ref{i4} $\Longrightarrow$ \ref{i1} Write $\G_n = \G_{\p_n}.$ We will
define a $\bT$ action on $\fA$ first by defining it on matrix units.
Here it will be convenient to regard $\bT$ as $\{ z \in \bC: |z| = 1
\}.$

Let $v$ be a matrix unit in $\fA;$ then (by compactness) $\h{v} =
\cup_{k=0}^{N} \G_k \cap \h{v},$ for some $N \in \bZ.$ Observe that
$\G_k$ is both open and closed in $\cP,$ so that, setting $\h{v}_k :=
\G_k \cap \h{v},$ each $v_k$ is a matrix unit or sum of matrix units
in $\fA.$ Now set \[ \al_z(v) = \sum_{k=0}^{N} e^{2\pi i kt}v_k,
\quad t \in \bT .\] $\al_t$ extends by linearity to (the dense
subalgebra of) all finite linear combinations of matrix units.  A
short calculation shows that $\al_z$ is isometric, so it extends to an
isometric map of $\fA.$ By the nested property of the $\G_n, $ the
automorphism property $\al_t(ab) = \al_t(a) \al_t(b) $ holds for $a,
b$ matrix units, hence for linear combinations of matrix units, and
finally for arbitrary $a, b \in \fA.$ One can verify directly, or use
\cite{ppw90} or \cite{scp92bk} to get that $\al_t$ is the restriction
of a star automorphism of C$^*(\fA).$ Finally, one notes that the
action $t \to \al_t$ is continous in the pointwise-norm topology;
i.e., for each $a \in \fA,$ the map $t \to \al_t(a)$ is norm
continuous.  \end{proof}

Recall from \cite{re94} that a continuous action $\al$ of $\bT$ on a
\cstar algebra $\fB$ is semi-saturated, if $\fB$ is generated as a
\cstar algebra by the fixed point algebra and the first spectral
subspace $\fB(1).$ (i.e., $\fB(1) = \{ b \in \fB:\ \al_t(b) = e^{2\pi
it}b .$)

\begin{theorem} \label{t:stdzanal} Let $\fA$ be a strongly maximal
TAF algebra with spectral triple $(X, \cP, \cR).$ The following conditions are
equivalent:
\begin{enumerate}
\item \label{i0'} There is a semi-saturated action $\al$ of $\bT$ on
C$^*(\fA)$
with fixed point algebra the diagonal $\fD,$ such that $\{ a \in
\fA:\ \al_t(a) = e^{2\pi int}a \} = (0),\ $ for $ n \geq 1.$
\item \label{i1'} There is a sequence $\fA(n), n \geq 0,$ of closed
linear subspaces of $\fA$ satisfying
\begin{enumerate}
\item \label{i1a} $\fA(0)$ is the diagonal $\fD;$ \item \label{i1b}
$\fA(n) \cap \fA(m) = (0),\ n \neq m;$
\item \label{i1c} $[\fA(n)
\cdot \fA(m)] = \fA(n+m),$ for all $n, m \geq 0;$ \item \label{i1d}
$\oplus_{n=0}^{\infty} \fA(n)$ is dense in $\fA.$ \end{enumerate}
\item \label{i2'} The counting cocycle $\h{d}$ is finite-valued and
continuous on $\cR,$ and $\h{d}^{-1}([0, \infty)) = \cP;$ \item
\label{i3'} $\fA$ is a standard TAF algebra; \item \label{i4'} There
is a partial homeomorphism $\p$ of $X$ such that the graph $\G_{\p}$
is open in $\cP \backslash \De(X),$ and \[ \cP \backslash \De(X) =
\cup_{n=1}^{\infty} \G_{\p^n} \] where $\p^n$ is the $n-$fold
composite of $\p.$ \end{enumerate}
\end{theorem}

\begin{proof}
\ref{i0'} $\Longrightarrow$ \ref{i1'}. With the spectral subspaces
$\fA(n)$ defined as in the proof of
\ref{t:zanal}, the condition that $\al$ be semi-saturated implies
that $\fA(n) = \text{closed span}\{a_1\cdots a_n:\ a_j \in \fA(1)\},\
n =2, 3, \dots.$ In particular, this
implies \ref{i1c}. Otherwise the proof is the same as \ref{i1}
$\Longrightarrow$ \ref{i2} of Theorem \ref{t:zanal}

\ref{i1'} $\Longrightarrow$ \ref{i0'}.  Define an action of $\bT$ on the
algebraic direct sum of the $\fA(n),\ n \geq 0,$ by \[ \al_t(\sum a_j)
= \sum e^{2\pi in_jt}a_j \quad \text{where } a_j \in \fA(n_j) \, .\]
$\al$ extends to an action of $\bT$ on $\oplus_{n\in \bZ} \fA(n),$
which is dense $\fA,$ hence to an action of $\fA.$ Define $\al_t$ on
$\fA(n)^* $ by $\al_t(a) = e^{-2\pi int}a, \ a \in \fA(n)^*,\ n > 0.$
This gives rise to an isometric action on $\fA^*,$ and hence a
star-action of $\bT$ on the norm-closure of $\fA + \fA^*,$ which is
$C^*(\fA).$

\ref{i1c} implies that the action is semi-saturated.

\ref{i1'} $\Longrightarrow$ \ref{i4'} We use the same notation as in
Theorem \ref{t:zanal}. Thus, $\h{\fA}(n) = \cup \{\h{v} :
v \text{ a matrix unit in } \fA(n) \}.$ By assumption, $\fA(1)\cap
\fA(0) = \fA(1) \cap \fD = (0),$ so that $\h{\fA}(1) \cap \De(X) =
\emptyset.$ So $\h{\fA}(1)$ is the graph $\G_{\p}$ of a partial
homeomorphism $\p$ of $X,$ and the graph $\G_{\p}$ is disjoint from
the diagonal set $\De(X).$ Also, repeated application of \ref{i1c}
implies that $\G_{\p^n} = \h{\fA}(n), \ n \geq 1.$ Now as in the
proof of \ref{i2} $\Longrightarrow$ \ref{i3} of Theorem
\ref{t:zanal}, \ $\cup_{n=0}^{\infty} \h{\fA}(n) = \cP,$ so that
$\cup_{n=1}^{\infty} \h{\fA}(n) = \cP \backslash \De(X).$ Hence,
$\cup_{n=1}^{\infty} \G_{\p^n} = \cP \backslash \De(X).$

\ref{i4'} $\Longrightarrow$ \ref{i1'} Let $\fA(n)$ be the closed
linear span of the set $\{ v \in \fA: v \text{ a matrix unit, } \h{v}
\subset \G_{\p^n} \}.$ We claim that \[\G_{\p^n} \cap \G_{\p^m} =
\emptyset, \text{ for } n \neq m. \quad \] Case i): $ n = 0 <
m.$ Of course if $m = 1,$\ it is true by assumption. If $\G_{\p^m}
\cap \De(x) \neq \emptyset$ for some $m > 1$, let $ (x, x) \in
\G_{\p^m} \cap \De(x).$ Then there is a $y \in X$ with $(x, y) \in
\G_{\p^{m-1}}$ and $(y, x) \in \G_{\p}.$ But then $(x, y)$ and $(y,
x)$ are both in $\cP,$ or $(x, y) \in \cP \cap \cP^{-1} = \De(X),$ so
that $y = x$ and $(x, x) \in \G_{\p},$ contrary to hypothesis. \\
Case ii) $ 0 < n < m.$ Suppose $(x, y) \in \G_{\p^n} \cap \G_{\p^m}.$
Then $\p^n(x) = y$ and $\p^m(x) = y.$ But $\p^m(x) =
\p^{m-n}\circ\p^{n}(x),$ so that $\p^{m-n}(y) = y.$ Since by Case i)
the graph of $\p^{m-n}$ is disjoint from the diagonal, this is
impossible. Thus the claim is established.

From the claim we have $\fA(n) \cap \fA(m) = (0)$ for $n \neq m.$
\ref{i1d} follows from the fact that $\oplus_{n=0}^{\infty} \fA(n)$
contains the algebra spanned by the matrix units of $\fA.$

\ref{i2'} $\Longleftrightarrow$ \ref{i3'} Follows from Proposition
2.8 and Theorem 2.9 of \cite{powa93}.

\ref{i2'} $\Longrightarrow$ \ref{i4'} Let $\G_{\p} = \h{d}^{-1}(1).$
This is open in $\cP$ and disjoint from $\De(X).$ It follows that
$\G_{\p^n} = \h{d}^{-1}(n),$ and hence \[ \cP = \h{d}^{-1}([0,
\infty)) = \cup_{n=0}^{\infty} \G_{\p^n}, \] so (with $\p^0 = id_X$),
$\cP \backslash\De(X) = \cup_{n=1}^{\infty} \G_{\p^n}.$

\ref{i4'} $\Longrightarrow$ \ref{i2'} Observe that the condition that
$\G_{\p}$ is open in $\cP \backslash \De(X)$ implies that the sets
$\G_{\p^n}$ are disjoint and open. Furthermore, the orbit (or
equivalence class) of a point $x \in X$ is given by $\cO(x) = \{
\p^n(x) : x \in \text{dom}(\p^n), n \in \bZ \}$ (where for $n$
negative, $\p^n(x)$ denotes the $n-$fold composite of $\p^{-1}$ at
$x.$) In particular, each orbit has the order type of a subset of
$\bZ,$ so that $\h{d}$ is finite on the groupoid $\cR: \ \h{d}(x, y)
= n$ if and only if $ y = \p^n(x).$ Thus the counting cocycle $\h{d}$
is finite and continuous on $\cR,$ and $\cP = \h{d}^{-1}([0,\infty)).$
\end{proof}

Let $X$ be a zero-dimensional compact space and $(X, \cP, \cR)$ a
spectral triple defined by a partial homeomorphism (respectively, a
nested sequence of partial homeomorphisms) on $X.$ Then Theorem
\ref{t:stdzanal} (respectively, \ref{t:zanal}) shows that $\cA(\cP)$ is
standard $\bZ$-analytic (respectively, $\bZ$-analytic) if and only if $\cR$
is an AF groupoid. In the next two sections we will give necessary
and sufficient conditions for $\cR$ to be AF.

\section{Ordered Bratteli diagrams and partial dynamical systems}
\label{s:pds}

Let $V$ and $W$ be two non-empty finite sets. An {\it ordered
diagram\/} from
$V$ to $W$ consists of a partially ordered set $E$ and surjective
maps $r: E\ra W$ and $s:E\ra V$ such that $e$ and $e'$ are comparable
iff $r(e) = r(e')$. Sometimes we just write $E$ for $(E,r,s)$. The
elements of $V$ and $W$ are the {\it vertices\/} and the elements of
$E$ are the {\it edges\/} of the diagram.

An {\it ordered
Bratteli diagram} $\ve$ consists of a vertex set \[
\cV = V_0 \cup V_1 \dots \text{ (disjoint union of finite sets)},
\]
where $V_0$ is a singleton, and
\[
\cE=\{ (E_n,r_n,s_n): n\ge 1\},
\]
where $(E_n,r_n,s_n)$ is an ordered diagram from $V_{n-1}$ to $V_n$.
If $e_i\in E_i$ with
$r_{i-1}(e_{i-1})=s_i(e_i)$ for all $i$, $m< i\le n$, then
$(e_{m+1},e_{m+2},\dots,e_n)$ is a {\it path\/} from $V_m$ to $V_n$.
Let $X=X\ve $ consist of all infinite sequences $(e_1,e_2,\cdots)$ of
edges with
$e_i\in E_i$ and $r_{i-1}(e_{i-1})=s_i(e_i)$ for all $i$. For
$x=(e_{n})_{n=1}^{\infty} \in X$ and $n\ge 1$, we will write $x(n) =
e_{n}$.

Suppose
$\ve$ is an ordered Bratteli diagram. For each path $p
=(e_{1},\,e_{2},\cdots,e_{n})$ from $V_{0}$ to $V_{n}$, let $C(p) =
\{(f_{1},\,f_{2},\cdots)\in X: f_{i}= e_{i}\text{ for all }1\le i\le
n\}$. We
give $X $ the smallest topology where each $C(p)$ is open. In this
topology,
each $C(p)$ is actually both closed and open (clopen).

\begin{proposition}
$X$ is a separable compact metrizable space. \end{proposition}

\begin{proof} Each $E_{n}$ is a discrete space, hence metrizable. Let
$Y = \prod_{n=1}^{\infty}E_{n}$, with the product topology.
Therefore, $Y$ is compact and metrizable and $X\subseteq Y$. The
topology of $X$ is equal to that inherited from $Y$. $X$ is separable
because it has a countable base $\{C(p): n\ge 1,\,p
=(e_{1},\,e_{2},\cdots,e_{n})\}$.
\end{proof}

Given an ordered Bratteli diagram $\ve$, let $X_{\max}$ to be the set
of
maximal paths, i.e., $X_{\max}=\{(e_i)\in X\ : e_i\text{ is maximal
in }E_i\text{ for all }i\}$. Similarly, define $X_{\min} = \{(e_i)\in
X\ : e_i\text{ is minimal in }E_i\text{ for all }i\}$. Since $X$ is
compact, it follows that these sets are always nonempty. Also it is
clear that the sets $\xmax,\ \xmin$ are closed. Now for every
$(e_i)\in X\ba \xmax$, let $k = \min\{i\ : e_i\text{ is not maximal
in }E_i\}$ and let $f_k$ be the successor of $e_k$ in $E_k$. For
$1\le i < k$, define $f_i$ so that $(f_1,\dots,f_{k-1})$ is the
unique minimal path (i.e., each $f_i$ is minimal in $E_i$)
from $V_0$ to $V_{k-1}$ such that
$r_{k-1}(f_{k-1}) = s_k(f_k)$.
Finally, let $f_n = e_n$ for $n>k$. Define a partial mapping $\p$ on
$X$
by $\p((e_i)) = (f_i)$. For
$x,\,y\in X$, we will write $x\le y$ if $\p^{n}(x) = y $ for some
$n\ge 0$. $\le$ is a partial ordering on $X$.

Let $x=(e_i)\in X$ and $n\ge 1$. Define $C_{n}(x) = C(p)$, where $p =
(e_{1},\cdots,\,e_{n})$. If $x=(e_i)\in X\setminus X_{\max}$, then
for every $n\ge 1$, there exists $m\ge n$ such that
$(e_{1},\,e_{2},\cdots,e_{m})$ is not maximal.  We have $\p(C_{k}(x))=
C_{k}(\p(x))$ for all $k\ge m$.  Since every open subset of
$X\setminus X_{\max}$ is a union of clopen subsets $C_{k}(x)$, $\p$ is
a partial homeomorphism.

\begin{definition} Let $\fB = (\cV, \cE, \geq)$ be an ordered
Bratteli diagram, and let $(X, \xmax, \xmin, \p)$ be the \pds
constructed above. The partial homeomorphism $\p$ is called a
Ver\v{s}ik transformation, and
$(X, \xmax, \xmin, \p)$ the Ver\v{s}ik \pds associated with $\fB.$
\end{definition}

For $x\in X$, let
\begin{eqnarray*}
\cO(x) &=& \{\p^{n}(x) : n\in \bZ, x \in \text{ dom } \p^n \}, \\
\cO^{+}(x) &=& \{\p^{n}(x) : n\ge 0, x \in \text{ dom } \p^n \}=
\{y\in X: y\ge x\},\\
\cO^{-}(x) &=& \{\p^{n}(x) : n\le 0, x \in \text{ dom } \p^n \}=
\{y\in X:y\le x\},
\end{eqnarray*}
and $w^{+}(x) $ and $w^{-}(x) $ be the accumulation points of
$\cO^{+}(x)$ and $\cO^{-}(x) $.

\begin{proposition} \label{p:orbit} Let $\ve $ be an ordered Bratteli
diagram and $X_{\max}$, $X_{\min}$, and $\p$ as defined above.
Suppose $U$ (respectively, $V$) is a clopen subset of $X$ containing
$X_{\min}$ (respectively, $X_{\max}$). Then we have \begin{eqnarray*}
i) &\cup_{n=0}^{\infty}\p^{n}(U)&=X\\
ii)&\cup_{n=0}^{\infty}\p^{-n}(V)&=X
\end{eqnarray*}
\end{proposition}

\begin{proof}
We prove only i), as the proof of ii) is similar. Let $x\in X$.  If
$\cO^{-}(x)$ is finite, then $\p^{-n}(x)\in X_{\min}$ for some $n\ge
0$. Therefore, $x\in \p^{n}(X_{\min})\subseteq \p^{n}(U)$. So, we may
assume that $\cO^{-}(x)$ is infinite.

We are going to show
that there exists a sequence $n_{1}> n_{2} >\cdots $ such that
$\p^{-n_{i}}(x)$ converges to some $x_{0}\in X_{\min}$.  Choose $m$ such that
$C_{m}(x_{0})\subseteq U$.  Then we have
$\p^{-n_{i}}(x)\in C_{m}(x_{0})\subseteq U$ for sufficiently large $i$, and hence $x\in
\p^{n_{i}}(C_{m}(x_{0}))\subseteq \p^{n_{i}}(U)$.

For each $k\ge 1$, let $(e_{1},\cdots, e_{k})$ be the unique minimal
path from
$V_{0}$ to $V_{k}$ such that $r(e_{k}) = r(x(k))$. Let $x^{k}\in X$ be
the point satisfying $x^{k}(i) =e_{i}$ for $1\le i \le k$ and
$x^{k}(i) = x(i) $ for $i>k$. Then $x^{k}= \p^{-n(k)}(x)$ for some
$n(k)>0$. Since $\cO^{-}(x)$ is infinite, $\{x^{k} :k\ge 1\} $ is
also infinite. Let $x_{0}$ be an accumulation point of $\{x^{k}:k\ge
1\} $. Then $x_{0}\in X_{\min}$ and there exists a subsequence
$n_{i}= n(k_{i})$ such that $\p^{-n_{i}}(x)$ converges to $x_{0}$.
\end{proof}

We will be studying the relationship of Ordered Bratteli diagrams and
partial dynamical systems. For the rest of this section, we will
mainly consider partial dynamical system satisfying the following
conditions:

\begin{definition}\label{d:bs}
	A Bratteli system is a quadruple $(X,\, \xmax,\, \xmin,\,\p)$ where
$\xmax$ and $ \xmin$ are closed subsets of the zero dimensional
compact space $X$, and $\p: X\setminus \xmax\to X \setminus \xmin$ is
a homeomorphism such that if $U$ (respectively, $V$) is a clopen
subset of $X$ containing $\xmin$ (respectively, $\xmax$), then
\begin{eqnarray*}
i)&\cup_{n=0}^{\infty}\p^{n}(U)&=X\\
ii)&\cup_{n=0}^{\infty}\p^{-n}(V)&=X
\end{eqnarray*}
Two Bratteli systems $(X^{1},\, X^{1}_{\max},\,
X^{1}_{\min},\,\p_{1})$ and
$(X^{2},\, X^{2}_{\max},\, X^{2}_{\min},\,\p_{2})$ are said to be
conjugate to each other if there exists a homeomorphism $h:X^{1}\to
X^{2}$ such that
$ h(X^{1}_{\max})= X^{2}_{\max}$, $h( X^{1}_{\min}) = X^{2}_{\min}$,
and $h\circ \p_{1}=\p_{2}\circ h$.
\end{definition}

\begin{remark} It will be shown that every Bratteli system is
conjugate to one arising from an ordered Bratteli diagram. This,
together with
Proposition \ref{p:orbit}, justifies the above definition.

Since the clopen sets form a base for the topology of $X,$ one could
just as well substitute `open' for `clopen' in the definition
\ref{d:bs}. However, mostly it will be convenient to work with clopen
sets.
\end{remark}

\begin{lemma}
Let $Y$ be a clopen set containing $\xmin$. Set $Z=\p^{-1}(Y) \cup
\xmax$. Then $Z$ is clopen. \end{lemma}

\begin{proof}
This follows from the fact that $\p|_{X\ba Z}$ is a homeomorphism
from $X\ba Z$ to $\p(X\ba Z) = X\ba Y.$

\end{proof}

Similarly, we have

\begin{lemma} Let $Z$ be a clopen set containing $\xmax$. Set
$Y=\p(Z) \cup \xmin$. Then $Y$ is clopen. \end{lemma}

We note that in the proof of the last two lemmas, the conditions (i),
(ii) in Definition \ref{d:bs} are not needed.

\begin{remark} Let $Z$ be a clopen set containing $\xmax$. Then $Y=\p(Z)
\cup \xmin $ is clopen and $Z=\p^{-1}(Y) \cup
\xmax$.   Similar result holds for $\xmin.$
\end{remark}

\begin{proposition} \label{p:bs}
Condition i) in the definition of Bratteli system is equivalent to
condition ii).
\end{proposition}
\begin{proof} We show i) $\Ra$ ii). By the remark, we may assume that $V =
\p^{-1}(Y)\cup \xmax$ for some clopen $Y\supseteq \xmin$. Suppose ii)
fails to hold.  Then there exists $x_{0} \in X$ such that for all $k \geq
0,\ \p^k(x_{0}) \notin V.$ In particular, $\p^k(x_{0}) \notin \xmax,$ so the
forward orbit $\{ \p^k(x_{0}):\ k \geq 0 \}$ is defined.  Since $V$ is
open, the closed orbit $\text{cl}\{ \p^k(x_{0}):\ k \geq 0 \} $ does not
intersect $V$, and hence $\om^+(x_{0})\cap V = \emptyset.$

Let $x \in \om^+(x_{0}).$ Thus there is a sequence $0 \leq k_1 < k_2 <
\cdots$ with $ x = \lim \p^{k_n}(x_{0}).$ By assumption i), $ x =
\p^m(y')$ for some $m \geq 0,\ y' \in Y.$ Thus, $y' = \lim_n \p^{k_n
- m}(y).$ Since $Y$ is open, there is an $n$ with $k_n - m > 0$ for
which $\p^{k_n -m}(x_{0}) \in Y.$ But then $\p^{k_n - m -1}(x_{0}) \in V.$
This is a contradiction, and the proof is complete.

The implication ii) $\Ra$ i) is analogous. \end{proof}

\begin{remark} \label{periodic} Note that Bratteli systems do not
admit periodic points; that is, there is no point $x \in X$ and
positive integer $n$ such that $x \in \textrm{ dom } \p^n$ and
$\p^n(x) = x.$
\end{remark}
\begin{proof} Assume to the contrary there is a periodic point $x$
with period $n.$ Then the orbit
$\cO(x) = \{ \p^j(x)\ : \ 0 \leq j < n \}$ is finite and evidently
disjoint from $\xmax.$ Thus there is a clopen neighborhood $Z$ of
$\xmax$ which is disjoint from $\cO(x).$ But by assumption ii) of
Bratteli systems there is a nonnegative integer $j$ with $x \in
\p^{-j}(Z);$ \ i.e., $\p^j(x) \in Z,$ a contradiction.
\end{proof}

\subsection{Bratteli diagrams from partial dynamical systems}

In this section we will show that any \pds satisfying (i), (ii)
(i.e., a Bratteli system)
is given by a Bratteli diagram. The proof follows Putnam's
construction in \cite{ifp89}.

Suppose $(X,\, X_{\max},\,
X_{\min},\,\p)$ is a Bratteli system. Let $Y$ be a clopen subset
containing $X_{\min}$ and $Z= \p^{-1}(Y) \cup \xmax$. Define
$\lambda:Y\to \bZ$ by \[
\lambda(y) =\min\{k\ge 0: \p^{k}(y) \in Z\}. \]
By condition (ii) of a Bratteli system, $\lambda(y)< \infty$ for all
$y\in Y$. By the
compactness of $X$, there is an $N\in\bZ^{+}$ such that
$\cup_{n=0}^{N}\p^{-n}(Z)= X$. Therefore, $\lambda(Y)$ only takes a
finite number of values, say, $J_{1}<J_{2}<\cdots < J_{m}$. Set
\[Y_{k}= \lambda^{-1}(J_{k})\subset Y,\quad \text{and}\quad Y(k,j) =
\p^{j}(Y_{k}),\quad \text{for}\quad j=0,\dots,J_{k}\, , \]
and $k=1,\dots, \, m$. Then we have
\begin{itemize}
\item[1.] $\cup^{m}_{k=1}Y(k,1) =\p(Y)$. \item[2.]$\p(Y(k,j)
)=Y(k,j+1)$, for $0\le j<J_{k}$. \item[3.]$\cup^{m}_{k=1}Y(k,J_{k})
=Z$.
\item[4.] $\cup^{m}_{k=1}\cup^{J_{k}}_{j=0}Y(k,j)=X$. \end{itemize}

It follows from the definition that the sets $Y(k,\,j),\,1\le k \le
m,\,0\le j\le J_{k}$ are disjoint.

1) and 2) are clear. To show 3), let $z\in Z$. By assumption (i),
$\cup_{n=0}^{\infty}\p^{n}(Y)=X$, so there is a smallest nonnegative
integer $j$ for which $z\in \p^{j}(Y)$. Set $y= \p^{-j}(z)$.

{\bf Claim} For $0\le \ell < j$, $\p^{\ell}(y) \notin Z$.

Suppose to the contrary that for some $\ell$, $0\le \ell < j$,
$\p^{\ell}(y)\in Z$. Note that $\p^{\ell}(y) \notin \xmax$, since
$\p^{\ell}(y) \in \textrm{dom }\ \p^{j-\ell}$. Thus, by definition of
$Z$,
$y_{1}=\p(\p^{\ell}(y))\in Y$, and $z=\p^{j-\ell-1}(y_{1})$. Since
$j-\ell-1\ge 0$, this contradicts our choice of $j$. This establishes
the claim.

Thus, $\lambda(y) =j$, so $j\in\{J_{1},\cdots,J_{k}\}$, and hence
$z\in\p^{J_{k}}(y)\in Y(k,\,J_{k})$ for some $k$.

The proof of 4) is similar. Let $x\in X$, and let $j$ be the
smalllest nonnegative integer with $x\in \p^{j}(Y)$. Set
$y=\p^{-j}(x)$. Then, as in the proof of 3), $\p^{\ell}(y)\notin Z$
for $0\le \ell < j$. If $y\in Y_{k}$, then $x\in Y(k,\,j)$, since
$0\le j\le J_{k}$.

Note: The sets $J(k,\,j)$, $1\le k \le m$, $0\le j\le J_{k}$ are open
and form a partition of $X$. Hence, each of them is both closed and
open (clopen).

\begin{definition}We will refer to such a partition as a
Kakutani-Rohlin partition, and
to the sets $\{Y(k,\,j) : 0\le j \le J_{k}\}$ as a tower.
\end{definition}

Let the \pds $(X, \xmax, \xmin, \p)$
be a Bratteli system, and let $\{Y_n\}_{n=1}^{\infty}$ be a nested
sequence of clopen sets containing $\xmin$ such that
$\cap_{n=1}^{\infty} Y_n = \xmin$. Then $\cap_{n=1}^{\infty}\left(
\p^{-1}\left(Y_n\right)\cup \xmax\right) = \xmax$. Let $\{\cP_n
\}_{n=1}^{\infty}$ be a nested sequence of finite clopen partitions
of $X$ such that $Y_n \in \cP_n,\ n = 1, 2, \dots,$ and $ \vee\cP_n$
is a base for the topology of $X.$

Inductively, construct sequences $\cQ_n,\ \cP'_n$ where $\cQ_n$ is a
finite clopen partition of $Y_n,$ \ $\cP'_n$ is a finite clopen
partition of $X$ as follows: set

\[ \cQ_1 = \bigvee_{j=0,\dots, J(1,k)}\bigvee_{k = 1,\dots, m(1)}
\{\p^{-j}[(Y_1(k,j)\cap P] \ : P \in \cP_1 \}. \]

So $\cQ_1$ is a partition of $Y_1,$ and each set $Y_1(k,0)$ is the
union of sets $Y_1(k,0)\cap Q$ as $Q$ runs through $\cQ_1.$ Index
these sets $Y_1(k,0,i),\ 1 \leq i \leq r(1,k).$ The sets
\[Y_1(k,j,i) = \p^j(Y_1(k,0,i))\]
$1 \leq k \leq m(1),\, 1 \leq i \leq r(1,k), 0 \leq j \leq J(1,k)$,
partition $X.$ Denote this partition $\cP_1'.$ Let $\cP_{0}'=\{X\}$.

Suppose now that $n > 1$ and $\cP_1', \dots, \cP_{n-1}'$ have been
defined
so that $\cP_l'$ is a refinement of $\cP_l$ and $\cP_{l-1}',$ for $ 1
\leq l \leq n-1.$ Set

\[ \cQ_n = \bigvee \left\{\p^{-j}(Y_n(k,j)\cap P)\ : \
\begin{array}{l}
P \in \cP_n \vee \cP_{n-1}', \\
1 \leq j \leq J(n,k), 1 \leq k \leq
m(n)
\end{array}
\right\} .\]

Thus the set $Y_n(k,0)$ is a union of sets $Y_n(k,0)\cap Q \ \ (Q \in
\cQ_n).$ Index these sets by $Y_n(k,0,i)\ , 1 \leq i \leq r(n,k)$.
Then, let $\cP_n'$ be the partition consists of the sets \[
Y_n(k,j,i) = \p^j(Y_n(k,0,i))\] $\ 1 \leq k \leq m(n),\ 1 \leq i \leq
r(n,k),\ 1 \leq j \leq J(n,k) $.

Note that since $\cP_n'$ is finer than $\cP_n,$ \ $\vee\cP_n'$
generates the topology of $X.$

\begin{theorem} \label{t:pds}
Let the \pds $(X, \xmax, \xmin, \p)$ be a Bratteli system. Then there
is a Bratteli diagram
$\fB = \fB(\cV, \cE ,\geq)$ such that the Ver\v{s}ik \pds $(X',
\xmax', \xmin', \p')$ associated to $\fB$ is conjugate to $(X, \xmin,
\xmax, \p).$
\end{theorem}

\begin{proof} The Bratteli
diagram $\fB$ to be constructed will have as its vertices the set \[
V_{n}=\{ Y_n(k,0,i)\ :\ 1 \leq i \leq r(n,k),\ 1 \leq k \leq m(n)
\}\quad\text{for }n\ge 1. \]

For convenience of notation, define $V_{0}=\{Y_{0}(1,0,1)\}=\{X\}$.
Let
$n\ge 1$. For each $ Y_n(k,0,i)$ and $0\le \ell\le J(n,k)$, if
$Y_n(k,\ell,i) \subseteq Y_{n-1}(k',0,i')$ for some $1\leq k' \leq
m(n-1)$, $1 \leq i' \leq r(n-1,k')$, we put an edge $e$ from
$Y_{n-1}(k',0,i')$ to $Y_n(k,0,i)$. We will indicate the
correspondence between $e$ and $\ell$ by $j(e)=\ell$. Define
$s(e)=Y_{n-1}(k',0,i')$ and $r(e)=Y_n(k,0,i)$. For two edges $e_{1}$,
$e_{2}$ with $r(e_{1})=r(e_{2})$, we put $e_{1}\le e_{2}$ if
$j(e_{1})\le j(e_{2})$. This defines an ordered Bratteli diagram
$\fB(\cV, \cE)$. Let $(X', \xmax', \xmin', \p')$ be the corresponding
Ver\v{s}ik \pds. We are going to show that $(X', \xmax', \xmin',
\p')$ is conjugate to $(X, \xmin, \xmax, \p).$

Since $Y_{n}\subseteq Y_{n-1}$, we have $j(e)=0$ for all minimal
edges $e$.
Suppose $e$ is a maximal edge from $Y_{n-1}(k',0,i')$ to
$Y_{n}(k,0,i)$. Let $\ell =j(e)$. Then we have
$Y_{n}(k,\ell,i)\subseteq Y_{n-1}(k',0,i')$ but $Y_{n}(k,j,i)\cap
Y_{n-1}=\emptyset$ for all $\ell <j\le J(n,k)$. Since for each $0\le
j \le J(n-1,k')$, $Y_{n-1}(k', j,i')\subseteq Y_{n-1}
$ or $Y_{n-1}(k', j,i')\cap Y_{n-1}=\emptyset $, we have $Y_{n-1}(k',
j-\ell,i')\cap Y_{n-1}=\emptyset $ for $\ell <j\le J(n,k)$. On the
other hand, $\p(Y_{n-1}(k', J(n,k)-\ell,i'))\cap Y_{n-1}\supset
\p(Y_{n}(k,J(n,k),i))\cap Y_{n}=\p(Y_{n}(k,J(n,k),i))\neq
\emptyset$. Hence, $J(n,k) -\ell = J(n-1,k')$, which implies that
$\ell = J(n,k) - J(n-1,k')$.

For $1\le m \le n$, let $e_{m}$ be an edge from
$Y_{m-1}(k_{m-1},0,i_{m-1})$ to $Y_{m}(k_{m},0,i_{m})$, then
$p=(e_{1},\,\cdots, e_{n})$ is a finite path in $\fB$ from $V_{0}$ to
$V_{n}$. For each $1\le m \le n$, let $j_{m}=j(e_{m})$. Define
$\Psi(p)= Y_{n}(k_{n}, s_{n}, i_{n})$ where
$s_{n}=\sum_{m=1}^{n}j_{m}$. Suppose $n>1$. Let $p' =
(e_{1},\,\cdots, e_{n-1})$ and $s_{n-1}=\sum_{m=1}^{n-1} j_{m}$. We
have
\begin{multline*}
\Psi(p) = Y_{n}(k_{n}, s_{n-1}+j_{n}, i_{n})
=\p^{s_{n-1}}(Y_{n}(k_{n}, j_{n}, i_{n}))\\
\subseteq \p^{s_{n-1}}(Y_{n-1}(k_{n-1}, 0, i_{n-1}))=\Psi(p').
\end{multline*}
Define
a map $\psi : X' \to X$ as follows. Let $x' \in X'$ be the infinite
path $x' = (e_1, e_2, \dots)$. Then let $\{\psi(x')\}=
\cap_{n=1}^{\infty}
\Psi((e_{1},\cdots,\,e_{n}))$. Clearly, $\psi$ is a homeomorphism
such that $\psi(X'_{\min})=\xmin$ and $\psi(X'_{\max})=\xmax$. Let
$x' = (e_1, e_2, \dots)\in X'$ such that $e_{\ell}$ is an edge from
$Y_{\ell-1}(k_{\ell-1},0,i_{\ell-1})$ to
$Y_{\ell}(k_{\ell},0,i_{\ell})$ and $j_{\ell}=j(e_{\ell})$ for all
$1\le \ell \le n$. If $x'$ is not maximal, let $m$ be the smallest
integer such that $e_{m}$ is not maximal. Let $j'_{m}$ be the
smallest integer $j> j_{m}$ such that $Y_{m}(k_{m}, j,
i_{m})\subseteq Y_{m-1}(k'_{m-1}, 0, i'_{m-1})$ for some $k'_{m-1}$,
$i'_{m-1}$. Since $Y_{m}(k_{m}, j, i_{m})\cap Y_{m-1}= \emptyset $
for all $j_{m}<j < j'_{m}$ and $Y_{m}(k_{m}, j'_{m}, i_{m})\subseteq
Y_{m-1}$, we have $j'_{m}-j_{m}=1+ J(m-1,k_{m-1})$. For $1\le \ell <
m$, $j_{\ell} = J(\ell,k_{\ell}) - J(\ell-1,k_{\ell}')$ because
$e_{\ell} $ is maximal. Therefore, $j'_{m} =1+
J(m-1,k_{m-1})+j_{m}=1+\sum_{\ell=1}^{m}j_{\ell}$. Let $f_{m}$ be the
edge from $Y_{m-1}(k'_{m-1}, 0, i'_{m-1})$ to $Y_{m}(k_{m}, 0,
i_{m})$ with $j(f_{m}) = j'_{m}$. Let $(f_{1},\cdots, f_{m-1})$ be
the unique minimal path from $V_{0}$ to $Y_{m-1}(k'_{m-1}, 0,
i'_{m-1})$. Then $\p'(x') = (f_{1},\dots, f_{m},
e_{m+1},e_{m+2},\dots)$. For $n>m$, let $s_{n} = \sum_{\ell
=m+1}^{n}j_{\ell}$. We have
\begin{multline*}
\{\psi(\p'(x'))\} =
\cap_{n=m+1}^{\infty}Y_{n}(k_{n},\sum_{\ell=1}^{m}j(f_{\ell})+s_{n},i_{n})\\
=\cap_{n=m+1}^{\infty}Y_{n}(k_{n},j'_{m}+s_{n},i_{n})\\
=\p(\cap_{n=m+1}^{\infty}Y_{n}(k_{n},\sum_{\ell=1}^{m}j_{\ell}+s_{n},i_{n}))
=\{\p(\psi(x'))\}.
\end{multline*}
It follows that $\psi\circ\p'=\p\circ\psi$.. \end{proof}

\begin{remark} It follows that for \emph{any} choice of nested clopen
sets $Y_n$ with intersection $\xmax$ the ordered Bratteli diagram
with corresponding Ver\v{s}ik transformation is conjugate to the
given \pds. Thus, the Ver\v{s}ik transformation is independent of the
nested sequence $\{ Y_n\}.$ \end{remark}

\begin{corollary} \label{c:afgpd}
Let $(X, \xmax, \xmin, \p)$ be a \pds. Then the following are equivalent:
\begin{enumerate}
\item $\cR = \{(x, \p^n(x)): x \in dom(\p^n),  n\in \bZ \}$ is an AF groupoid; 
\item the \pds $(X, \xmax, \xmin, \p)$ is conjugate to Ver\v{s}ik map on a Bratteli compactum;
\item for any clopen subset $U$ containing $\xmin, \ \cup_{n=0}^{\infty} \p^n(U) = X;$
\item for any clopen subset $V$ containing $\xmax, \ \cup_{n=0}^{\infty} \p^{-n}(V) = X.$
\end{enumerate}
\end{corollary}

For the next result we will need the notion of equivalence for
ordered Bratteli diagrams from \cite{powa93} (Definitions 3.4, 3.5
and 3.6); for the reader's convenience, we recall the definition here.

\begin{definition} Let $V,\ W$ be (finite) vertex sets. Two ordered
diagrams $(V, E, r, s), \ (V, E', r', s')$ are \emph{order
equivalent} if there is an order-preserving bijection $\Phi:\ E \to
E'$ such that

\[ r(e) = r'(\Phi(e)) \quad \text{and} \quad s(e) = s'(\Phi(e))\, .\]

Now, let $(\cV, \cE),\ (\cW, \cF)$ be two ordered Bratteli diagrams.
The Bratteli diagrams are said to be order equivalent if there exist
strictly maps $g, h: \bN \to \bN$ with $g(0) = h(0)$ and ordered
diagrams $E_n'$ from $V_n$ to $W_{g(n)}$ and $F_n'$ from $W_n$ to
$V_{h(n)}$ such that

\[ F_{g(n)}' \circ E_n' \text{ is order equivalent to }
E_{h(g(n))}\circ \cdots \circ E_{n+1} \]

and

\[ E_{h(n)}'\circ F_n' \text{ is order equivalent to } F_{g(h(n))}
\circ \cdots \circ F_{n+1}\ . \]
\end{definition}

\begin{remark} A contraction of an ordered Bratteli diagram $\fB =
(\cV, \cE)$ is another ordered Bratteli diagram $(\cV', \cE')$
together with a subsequence $\{n_k\}$ of the positive integers such
that $V_k' = V_{n_k}$ and the edge set $E_k'$ consists of all paths
from $V_{n_{k-1}}$ to $V_{n_{k}},$ ordered lexicographically.

Another definition of order equivalence of ordered Bratteli diagrams
is as follows: say $\fB = (\cV, \cE), \
\fB' = (\cV', \cE')$ are order equivalent if there is a third ordered
Bratteli diagram $\fB'' = (\cV'', \cE'')$ such that the contraction
of $\fB''$ to the even vertices is a contraction of $\fB,$ and the
contraction of $\fB''$ to the odd vertices is a contraction of $\fB'.$

It is easy to see that this definition is equivalent to the one we
have given above.
\end{remark}

While we have been viewing ordered Bratteli diagrams as a partial
dynamical systems, one can just as well view them as defining a
standard $\bZ$-analytic algebra by viewing the edges as (ordered)
standard embeddings, as in \cite{powa93}.  Now it is known that two
semicrossed products are isomorphic iff their actions
(i.e.,homeomorphisms) are conjugate (\cite{jp84}, \cite{haho88},
\cite{scp92a}). Exel's generalized notion of crossed product by a
partial action (\cite{re94}, \cite{re95}) allows one to consider
standard $\bZ$-analytic algebras as semicrossed products.  In this
sense, the following theorem extends the earlier results on
semicrossed products.

\begin{theorem} \label{t:conj}
Let $(\cV, \cE)$ (resp. $(\cV', \cE')$) be an ordered Bratteli
diagram, $(X, \xmax, \xmin, \p)$
(resp., $(X', \xmax', \xmin', \p')$) the Ver\v{s}ik \pds constructed
from the diagram, and $\fA = \fA(\cV, \cE)$ (resp., $\fA' =
\fA'(\cV', \cE')$) the standard $\bZ$-analytic TAF algebra defined by
the diagram $(\cV, \cE)$ (resp., $(\cV', \cE')$). Then the following
are equivalent: \begin{enumerate}
\item \label{a} The Bratteli diagrams $(\cV, \cE),\ (\cV', \cE')$ are
order equivalent;
\item \label{b}The Ver\v{s}ik partial dynamical systems $(X, \xmax,
\xmin, \p)$,
{}\newline $ (X', \xmax', \xmin', \p')$ are conjugate; \item
\label{c}$\fA$ is isometrically isomorphic to $\fA'.$ \end{enumerate}
\end{theorem}

\begin{proof} \ref{a} $\Longleftrightarrow$ \ref{c} was proved in
\cite{powa93}, Theorem 3.7.

\ref{b} $\Longrightarrow$ \ref{c} To begin, we need the fact that
if a strongly maximal TAF algebra is the inductive limit of a system
$(\fA_n, \sigma_n),$ then the spectrum $\cP$ of $\fA$ is the
projective limit of the spectra of $\fA_n.$ It follows in this case
that the spectrum $\cP$ of $\fA $ is $\{ (x, y) \in X\times X: \ y =
\p^n(x) , \text{ for some } n \geq 0 \}.$

Suppose \ref{b}. is satisfied, and let $\psi : X \to X'$ be a
homeomorphism which induces a conjugacy of the two partial dynamical
systems. Define
$\Psi : \cP \to \cP'$,\[ \Psi(x, \p^j(x)) = (\psi(x), \psi(\p^j(x)) =
(\psi(x), {\p'}^j(\psi(x)),\]
$x \in X\ba \xmax\, $. Then $\Psi$ is a
semigroupoid isomorphism, and so by \cite{scp92bk}, $\fA,\ \fA'$ are
isometrically isomorphic.

\ref{c} $\Longrightarrow$ \ref{b} Conversely, suppose there is a
semigroupoid isomorphism $\Psi : \cP \to \cP'.$ Setting $\Delta(X) =
\{ (x, x):\ x \in X \},$ note that $\Psi(\Delta(X)) = \Delta(X').$
Indeed, if $\Psi(x, x) = (x', y')$,  then
%\begin{multline*}
%\Psi(x, x) = (x', y'), \text{ then }
$$(x', y') = \Psi((x, x)\circ(x,
x))
= \Psi(x, x)\circ\Psi(x, x) = (x', y')\circ(x', y'),$$
%\end{multline*}
forcing $y' = x'.$ Thus there is a homeomorphism $\psi: X \to X'$ so
that $\Psi(x, x) = (\psi(x), \psi(x))\, ,\ x \in X.$ Also,
\begin{gather*}
\Psi(x, \p(x)) = \Psi(x, x)\circ\Psi(x, \p(x))\circ \Psi(\p(x),
\p(x)) \\
= (\psi(x), \psi(x))\circ \Psi(x, \p(x)) \circ(\psi(\p(x)),
\psi(\p(x)))\, ,
\end{gather*}
from which $\Psi(x, \p(x)) = (\psi(x), \psi(\p(x))).$ Since
$(\psi(x), \psi(\p(x)) \in \cP'$, $ \psi(\p(x)) = {\p'}^j(\psi(x))$ for
some $ j \geq 0.$ Now $j=0$ is impossible, as $\Psi^{-1} $ maps
$\Delta(X')$ into $\Delta(X)$.
Now if $j \geq 2,$ then
\[ (\psi(x), {\p'}^j(\psi(x))) = (\psi(x), \p'(\psi(x)))\circ
(\p'(\psi(x)), {\p'}^j(\psi(x))) \]
is the composition of two elements of $\cP' \ba \Delta(X').$ Applying
$\Psi^{-1},$ we have that $(x, \p(x))$ is the composition of two
elements of $\cP \ba \Delta(x).$ But this is impossible, as $\p(x)$
is the immediate successor of $x.$ Thus,
\[ \Psi(x, \p(x)) = (\psi(x), \p'(\psi(x)))\, ,\text{ and hence }
\p'(\psi(x)) = \psi(\p(x))\, . \]
In other words, $\psi$ is a conjugacy of the two partial dynamical
systems.
\end{proof}

\section{Nested Sequences and AF-algebras} \label{s:nested} In this
section, we are
going to characterize nested sequences of partial homeomorphism that
will give rise to AF-algebras.

Let $\fA$ be a $\bZ$-analytic algebra and $\fB =C^{*}(\fA)$. Choose a
sequence
$\{ \fB_n \}$ of finite dimensional \cstar algebras of $\fB$, say
$\fB_n = \bigoplus _{k=1}^{r(n)} M_{m(n, k)}$, and a set of matrix
units $\{ e_{ij}^{(nk)} \}$ for $\cup _{n=1}^{\infty} \fB_n$ such that
\[ \fA_n:= \fA\cap \fB_n = \bigoplus _{k=1}^{r(n)} T_{m(n, k)}, \]
where $T_{m(n, k)}$ is the upper triangular subalgebra of $ M_{m(n,
k)},$ and $\fA = \dirlim \fA_n$ and $\fB = \dirlim \fB_n$. Denote
$ e_{ii}^{(nk)}$ by $ e_{i}^{(nk)} $.
Let $\left( X,\cP,\cR\right) $ be the associated spectral triple.
Then there
exists an integer valued cocycle $d$ on $\cR$ which induces a nested
sequence $\left\{\p_{n}\right\}_{n\in \bZ} $ of partial homeomorphism
on $X$ such that
\[ \cP = \sqcup_{n=0}^{\infty} \G_{\p_n} .\]
For $n\in \bZ$, let $\xmax^{n}=\xmin^{-n}=X\setminus \dom \p_{n}$,
$\xmax = \cap_{n=1}^{\infty}\xmax^{n}$ and $\xmin =
\cap_{n=1}^{\infty}\xmin^{n}$.

We have
\[\xmin = \cap_{n}\left(\cup_{k}\hat
e_{1}^{(nk)}\right) \text{ and } \xmax = \cap_{n}\left(\cup_{k}\hat
e_{m(n,k)}^{(nk)}\right) \]

Suppose $U$ and $V $ are clopen subsets of $X$ containing $\xmin$ and
$\xmax$ respectively. Then there exists $N$ such that $\cup_{k}\hat
e_{1}^{(Nk)}\subseteq U$ and $\cup_{k}\hat e_{m(N,k)}^{(Nk)}\subseteq
V$. For $n\in\bZ$, let $\ti{\p}_n$ be the restriction of $\p_n$ to
$\cD_n= \{x \in X: \left(x,\p_n(x)\right)\in \hat{e}_{ij}^{(Nk)}
\text{ for some }i,\,j\}$. Then each $\cD_n$ is clopen. Let

\[M= \max\{d(x,y):(x,y)\in \hat e_{ij}^{(Nk)} 1\le i,\, j \le
m(N,k),\,1\le k \le r(N)\}. \]
If $q> M$ then the domain of $ \tilde
\p_q$ is empty. This proves the necessity of the conditions in
Theorem \ref{t:AFnest}. We first prove the sufficiency of the
conditions for a special case of the theorem.

\begin{lemma} \label{l:Phi} Let
$\cN= \left\{\p_n\right\}_{n\in \bZ}$ be a nest of partial
homeomorphism on compact zero-dimensional space $X$ such that
\begin{enumerate}
\item $\xmin$ and $\xmax$ are clopen.
\item For each $n\ge 1 $, $\dom \p_n $ and $\dom \p_{-n}$ are clopen.
\item There exists $M\ge 1$ such that $\dom \p_n= \emptyset$ for
$|n|\ge M$
\end{enumerate}
Then the groupoid defined by $\cN$ is AF. \end{lemma}
\begin{proof} Define $\Phi:X\setminus \xmax \to X\setminus \xmin $ by
\[
\Phi(x) = \p_k(x)\text{ where }k= \min\left\{n\ge
1:x\in\dom\p_n\right\}.
\]
We are going to prove that
\begin{enumerate}
\item $\Phi $ is bijective and \[
\Phi^{-1}(y) = \p_{-k}(y) \text{ where }k= \min\left\{n\ge
1:y\in\dom\p_{-n}\right\}.
\]
\item $\Phi$ is a homeomorphism. \item
$\cup_{n=0}^{\infty}\Phi^n\left(\xmin\right)= X=
\cup_{n=0}^{\infty}\Phi
^{-
n}\left(\xmax\right)$.
\end{enumerate}

{\bfseries Proof of 1.} Suppose $x,\,y\in \xmaxc$ and $\Phi(x)=
\Phi(y)$. Let
\[r=\min\left\{n\ge 1:x\in\dom\p_n\right\}\text{ and } s=
\min\left\{n\ge 1:y\in\dom\p_n\right\}. \]
So $\Phi(x)= \p_r(x)= \p_s(y)= \Phi(y)$.  If $r> s$, then we have $y=
\p_{- s}\circ \p_r(x)=\p_{r-s}(x)$.  Therefore, $x\in\dom \p_{r-s}$
and $1\le r-s < r$, a contradiction.  Therefore, $r\le s$.  Similarly,
$r\ge s$ and consequently, $r= s$.  This proves that $\Phi$ is one to
one.

Let $y\in \xminc$ and $k= \min\left\{n\ge 1:y\in\dom\p_{-n}\right\}$.
Then $y= \p_k(x)$ for some $x\in \xmaxc$. Suppose $n\ge 1$ and $x \in
\dom\p_n$. Then $y = \p_k(x)= \p_{k-n}\circ\p_n(x)\in \dom\p_{-(k-
n)}$ and $k-n< k$. Therefore, $k-n\le 0$. Hence, $n\ge k $ and $y=
\Phi(x)$.

{\bfseries Proof of 2.} Let $C$ be a closed subset of $\xmaxc$. We
will show that
$\Phi(C)$ is closed. Suppose $x_n\in C$ and $y=
\underset{n\to\infty}\lim\Phi(x_n)$. By choosing a subsequence if
necessary, we may assume that for some fixed $k\ge 1$, $\Phi(x_n)=
\p_k(x_n)$ for all $n$ and $x_0=\underset{n\to\infty}\lim x_n\in C\cap
\dom\p_k$.  Therefore, $y= \p_k(x_0)$.  Suppose $m\ge 1$ and $x_{0}\in
\dom\p_m$.  Then $x_n \in \dom\p_m$ for sufficiently large $n$.
Hence, $k\le m$ and $y=\p_k(x_0)=\Phi(x_0)\in \Phi(C)$.  Therefore,
$\Phi(C)$ is closed.  Since $\Phi$ is bijective, $\Phi(O)$ is open for
every open set $O$ in $\xmaxc.$ The same proof shows that
$\Phi^{-1}(O)$ is open for any open set $O$ in $\xminc.$ Consequently,
$\Phi$ is a homeomorphism.

{\bfseries Proof of 3.} Since $\dom\Phi^M=\emptyset$, we have
$$X=
X\setminus
\dom\Phi^M= \cup_{n=0}^{M-1}\Phi^{-1}(\xmax).$$ Hence,
$X=
\cup_{n=0}^{\infty}\Phi^{-n}\left(\xmax\right)$. Similarly,
$X=\cup_{n=0}^{\infty}\Phi^n\left(\xmin\right)$.

It follows that $\left(X,\,\xmax,\,\xmin,\,\Phi\right)$ is a Bratteli
system. By
construction, this system generates the same groupoid as that defined
by
the nest $\cN$.

\end{proof}

\begin{theorem} \label{t:AFnest}Let $\cN =\left\{\p_n\right\}_{n\in
\bZ}$ be a nest of partial homeomorphism on compact zero-dimensional
space $X$ and $\cR$ the groupoid defined by $\cN$. Then $\cR$ is AF
if and only if the following is satisfied:

Given clopen subsets $U$ and $V $ of $X$ containing $\xmin$ and
$\xmax$ respectively, there exist clopen subset $Y$ and $Z$ with
$\xmin \subseteq Y\subseteq U$ and $\xmax \subseteq Z\subseteq V$
such that
for $n\ge 1$, if $\ti{\p}_n$ is the restriction of $\p_n$ to
$\cD_{n}=\left\{x\in X\setminus Z:\p_n(x)\in X\setminus Y\right\}$ , $\tilde
\p_{- n}=\tilde \p_n^{-1}$, and $\p_0=id_X,$ then the system
$\{\ti{\p_n}\}$ satisfies the conditions of Lemma \ref{l:Phi}.
\end{theorem}

\begin{proof} Suppose $\cN$ is a nest satisfying the above
conditions. Then we can choose a sequence of clopen subsets $U^k$ and
$V^k$ such that $\cap_k U^k = \xmin$ and $\cap_kV^k = \xmax$. For
each $k$, applying the condition to $U=U^k$ and $V= V^k$, we have
clopen sets=
$Y^k $
and $Z^k$ with $\ti{\p}_n^{(k)}$ and $\dom \ti{\p}_n^{(k)}$ defined
accordingly.
Let $\cR^k $ be the groupoid defined by the nest
$\left\{\ti{\p}_n^{(k)}\right\}$. Since $\cR_k\subseteq \cR_{k+1}$ and
$\cR= \cup_{k=1}\cR_k$, it suffices to prove that each $\cR_k$ is AF.
Since the nest $\left\{\tilde\p_n^{(k)}\right\}$ satisfies the
conditions in Lemma \ref{l:Phi}, the result follows.
\end{proof}

\begin{example}\label{e:dh} This is an example of a nested sequence
of partial homeomorphisms satisfying the conditions in Theorem
\ref{t:AFnest} but the counting cocycle on the associated groupoid is
not continuous. Therefore, the corresponding TAF algebra is
$\bZ$-analytic but not standard $\bZ$-analytic. This TAF algebra is
isomorphic to the example given by Donsig and Hopenwesser
\cite{ppw94}.

Suppose $X
= \prod _{n=1}^{\infty} \{ 0, 1 \}$. Define a homeomorphism $\p:X\to
X$ by $\p((x_{n}))=( y_{n})$, where
\[
y_{n}=\left\{
\begin{array}{ll}
0&\ \ \text{ if }x_{i}=1 \text{ for all }1\le i\le n \\ 1&\ \ \text{
if }x_{i}=1 \text{ for all }1\le i\le n-1\text{ and }x_{n}= 0 \\
x_{n}&\ \ \text{ if }x_{i}=0 \text{ for some }1\le i\le n-1
\end{array}
\right.
\]
($\p$ is usually referred to as the odometer map.)

For $x\in X$, let $\xmax = \{x\}$ and $\xmin=\{\p(x)\}$. Restricting
$\p$ to $\xmaxc$, we have a \pds $\left(X,\,\xmax,\,\xmin,
\p\right)$, which is a Bratteli system. The conjugacy class of this
system is independent of the choice of $x$. Indeed, viewing $X$ as a
solenoidal group, given any two points $x, x'$
there is a homeomorphism $h$ of $X$ mapping $x$ to $x'$ which
commutes with $\p,$ namely $h(y) = y +x - x'.$ It follows that the
system with $x$ as the maximal point is conjugate to the system with
$x'$ as the maximal point. The corresponding TAF
algebra is standard $\bZ$-analytic with $\fB_{n}=M_{2^{n}}$ and
$\fA_{n}=T_{2^{n}}$.

Let $\mathbf{x_{0}}=(0,\,0,\,\cdots) $ and
$\mathbf{x_{n}}=\p^{n}\left(\mathbf{x_{0}}\right)$ for $n\in \bZ$. Define
$$
\begin{array}{rcl}
\xmax^{0}&=&\emptyset\\
\xmax^{1}&=&\left\{\mathbf{x_{k}}: k=-2,-1,0\right\}\\
%\xmin^{1}&=&\left\{\mathbf{x_{k}}:k=-1,0,1\right\}\\
\xmax^{2}&=&\left\{\mathbf{x_{k}}:k=-3,0\right\}\\
%\xmin^{2}&=&\left\{\mathbf{x_{k}}:k=-1,2\right\}\\
%&&\text{and\ for }n> 2, \\
\xmax^{n}&=&\left\{\mathbf{x_{k}}:-n-1\le k \le 0,\, k\ne
-n,-1\right\}\ \text{ for }n> 2,\\
%\xmin^{n}&=&\left\{\mathbf{x_{k}}:k=-1,1,2\cdots, n-2,n\right\}\\
\xmax^{-n} &=&\p^{n}(\xmax^{n})\ \text{ for }n \ge 1.\\
\end{array}
$$ Let
$\p_{0}= id_{X}$. For $n\ne 0$, let $\p_{n}= \p^{n}|_{X\setminus
\xmax^{n}}$. Then $\cN=\left\{ \p_n
\right\}_{n\in \bZ}$ is a nested sequence of partial homeomorphism on
$X$ with $\xmax = \{\mathbf x_{0}\}$ and $\xmin=\{\mathbf{x_{-1}}\}$. We
are going to
show that
\begin{enumerate}
\item $\cN$ satisfies the conditions in Theorem \ref{t:AFnest}. \item
The counting cocycle $\h d$ is not continuous. \end{enumerate}

{\bfseries Proof of 1.} Given $k\ge 1$ and $c_{i}\in \{0,\,1\}$ for
$i=1,\,\dots,\,k$, define the $k-$cylinder sets \[
\left[c_{1},\dots,\,c_{k}\right]_{k}= \left\{\mathbf x =\{x_{n}\}\in X:
x_{i}=c_{i}\ \text{ for }1\le i \le k\right\}. \]
The cylinder sets are clopen and form a basis of the topology of $X$.
Given any clopen sets $U$ and $V$ with $\xmin\subseteq U$ and
$\xmax\subseteq V$, there exists $k\ge 1$ such that
$[1,\,\dots,\,1]_{k}\subseteq U$ and
%$[\underbrace{0,\,\dots,\,0}_{\text k-terms}]\subseteq V$.
$[0,\,\dots,\,0]_{k}\subseteq V$. Let $\cC_{k}$ be the collection of
all $k-$cylinder sets. Then
\[
\cD_{n}=\left\{
\begin{array}{ll}
\cup\left\{C\in \cC_{k}:C\cap \xmax^{(n)}= \emptyset\right\}
&\quad\text{
for } n \ne 0\\
X&\quad \text{ for } n=0\\
\end{array}
\right.
\]
It is straightforward to check that condition in Theorem
\ref{t:AFnest} is satisfied.

{\bfseries Proof of 2} Let $\mathbf{x^{n}}=\left
\{x^{n}_{i}\right\}_{i\ge 1} $,
where $x^{n}_{i}= 1$ for $1\le i\le n$ and $x^{n}_{i}=0$ otherwise.
Then each $\mathbf{x^{n}}$ is in $\dom\p_{1}\circ \p_{1}$ but
$\mathbf{x^{n}}$ converges to $\mathbf{x_{-1}}\in \dom \p_{2}\setminus
(\dom\p_{1}\circ \p_{1})$. Therefore, $\h
d(\mathbf{x^{n}},\p_{2}(\mathbf{x^{n}}))=2$ but $\h
d(\mathbf{x_{-1}},\p_{2}(\mathbf{x_{-1}}))=1$. Hence, $\h d $ is not
continuous.
\end{example}

Let
$\cN=\left\{\p_n\right\}_{n\in \bZ}$ be a nest of partial
homeomorphism satisfying the condition in Theorem \ref{t:AFnest}.
Define
$\Phi:X\setminus \xmax \to X\setminus \xmin $ by \[
\Phi(x) = \p_k(x)\text{ where }k=\min\left\{n\ge
1:x\in\dom\p_n\right\}
\]
(as defined in Lemma \ref{l:Phi}).
In Example \ref{e:dh}, we have $\mathbf{x^{n}}\to \mathbf{x_{-1}}$
but $\Phi(\mathbf{x^{n}})\to \mathbf{x_{0}}\ne \mathbf{x_{1}}=
\Phi(x_{-1})$. Therefore, $\Phi$ may not be continuous. However,
$\Phi$ turns out to be very useful in the study of the isomorphism of
the associated $\bZ$-analytic algebra. The proof of the following
theorem is similar to the proof for the equivalence of conditions 2
and 3 in Theorem \ref{t:conj}.

\begin{theorem} \label{t:nconj}
Suppose $\cN=\left\{\p_n\right\}_{n\in \bZ}$ and
$\cN'=\left\{\p'_n\right\}_{n\in \bZ}$ are two nested sequence of
partial homeomorphisms on $X$ and $X'$ satisfying the conditions in
Theorem \ref{t:AFnest}.  Then the associated $\bZ$-analytic algebras
$\fA$ and $\fA'$ are isometrically isomorphic if and only if there
exists a homeomorphism $\psi:X\to X'$ such that $\psi(\xmax)=\xmax'$,
$\psi(\xmin)=\xmin')$ and $\psi\circ \Phi =\Phi'\circ \psi$.
\end{theorem}

\subsection{Semi-saturating a nested sequence}
%nonsemi-saturated action}
\label{ss:sat} In Example \ref{e:dh} there is a single partial
homeomrophism $\p$ of $X$ (actually in this case a homeomorphism) such
that each $\p_n$ in the nested sequence arises as a restriction of the
$n-$fold composition $\p^n.$ A \emph{semi-saturation} of a nested
sequence $\{ \p_n \}$ of partial homeomorphisms is a partial
homeomorphism $\p$ on $X$ satisfying $ dom(\p_n)\subset dom(\p^n),$
and for all $x \in dom(\p_n), \ \p_n(x) = \p^n(x), \ n = 1, 2, \dots.$

We note that in the above example, $\dom \p_{1}$ is dense in
$\cup_{n=1}^{\infty}\dom \p_{n} $. For
such systems we have the following result.

\begin{proposition} \label{p:semisat}
Let $\{\p_{n}\}_{n=0}^{\infty}$ be a nested sequence of partial
homeomorphism
on $X$ such that $\dom \p_{1}$ is dense in $\cup_{n=1}^{\infty}\dom
\p_{n}$. Then the following
conditions are equivalent:
\begin{itemize}
\item[1.] There exists a homeomorphism $\p$ on
$\cup_{n=1}^{\infty}\dom \p_{n} $
such that $\p_{n}=\p^{n}\big |_{\dom \p_{n}}$ for all $n\ge 1$.
\item[2.]
\begin{itemize}
\item[a] For every sequence $\{\mathbf {x^{(k)}}\}$ in $\dom \p_{1}$
such
that $ \underset{k\to
\infty}\lim \mathbf {x^{(k)}}\in \cup_{n=1}^{\infty}\dom \p_{n}$, we
have that $
\underset{k\to \infty}\lim \p_{1}(\mathbf{x^{(k)}})
%\in \cup_{n=1}^{\infty}\mbox{ range } \p_{n}
$ exists; and \item[b] for every sequence $\{\mathbf {y^{(k)}}\}$ in $
\mbox{ range } \p_{1}$ such that $ \underset{k\to \infty}\lim \mathbf
{y^{(k)}}
\in \cup_{n=1}^{\infty} \mbox{ range } \p_{n}$, we have that
$ \underset{k\to \infty}\lim \p_{1}^{-1}(\mathbf {y^{(k)}})
%\in \cup_{n=1}^{\infty}\dom \p_{n}
$ exists.  \end{itemize}
\end{itemize}
\end{proposition}

\begin{proof}
Clearly, 2 follows from 1. Suppose 2 holds. Let $\mathbf{ x} \in
\cup_{n=2}^{\infty}\dom \p_{n}$ and $\mathbf{x^{(k)}}\in \dom \p_{1}$
such that
$ \underset{k\to \infty}\lim \mathbf{x^{(k)}}=\mathbf{x}$. Then $\y=
\underset{k\to \infty}\lim \p_{1}(\mathbf{x^{(k)}})$ exists and the limit
is independent of the choice of $\{\mathbf{x^{(k)}}\}$.  Therefore, we
can extend $\p_{1 }$ to a continuous map $\p$ on
$\cup_{n=1}^{\infty}\dom \p_{n}$.  It follows from b that $\p$ is a
homeomorphism on $\cup_{n=1}^{\infty}\dom \p_{n} $.  \end{proof}

\begin{example} \label{e:nonsemisat}
Let $X= \{(\frac 1n,\,i):i=1,2{\rm \ and\ }n\ge 1\}\cup
\{(0,1),(0,2)\}\subset \bR^{2}$, $X_{1}=X\setminus \{(0,1),(0,2)\}.$
Define a nested sequence
$\{\p_{n}\}_{n=0}^{\infty}$ of partial homeomorphism on $X$ by:
$\p_{0}={\rm id}_{X}$; $\p_{1}((\frac
1{2k-1},\,i))=(\frac 1{2k},\,i))$ for $i=1,\,2$ , $\p_{1}((\frac
1{2k},\,1))=(\frac 1{2k-1},\,2))$, and $\p_{1}((\frac
1{2k},\,2))=(\frac 1{2k+1},\,1))$ for all $k\ge 1$;
$\p_{2}((0,\,1))=(0,\,2)$ and $\p_{2}(\mathbf{x})=\p_{1}^{2}(\x)$ for
all $\x\in X_{1}$ and for $k\ge 2$, $\p_{2k-1}=\p_{1}^{2k-1}$ and
$\p_{2k}=\p_{2}^{k}$. Let $\mathbf{x^{k}}= (\frac 1{2k-1},\,1) $ and
$\mathbf{y^{k}}= (\frac 1{2k},\,1) $. Then $ \underset{k\to
\infty}\lim \mathbf{x^{(k)}}= \underset{k\to \infty}\lim \mathbf
{y^{(k)}}= (0,\,1)$ but $ \underset{k\to \infty}\lim \p_{1}(\mathbf
{x^{(k)}})= (0,\,1)$ and $ \underset{k\to \infty}\lim \p_{1}(\mathbf
{y^{(k)}})= (0,\,2)$. So there exists no partial homeomorphism $\p $
on $X$ such that $\p_{n}=\p^{n}\big |_{X_{n}}$ for all $n$.

This shows that the nested sequence $\{ \p_n \}_{n=1}^{\infty}$ does
not admit a semi-saturation.
Let $(X,\,\cP_{1},\,\cR_{1})$ be the spectral triple associated with
$\{\p_{n}\}_{n=0}^{\infty}$. We want to show that $(X, \cP_1, \cR_1)$
is the spectral triple of a
$\bZ$-analytic algebra. For this it is enough to know that $\cR_1$ is
an AF-groupoid.

Define $\xmax =\{(0,2)\}$, $\xmin = \{(1,\,1)\}$.
and $\p:X\setminus \xmax \to X\setminus \xmin$ by $\p((0,\,1))=(0,\,2),\,\p((\frac
1{n},\,1))=(\frac 1{n},\,2))$ and $\p((\frac 1{n},\,2))=(\frac
1{n+1},\,1))$ for all $n\ge 1$. Then $(X,\,\xmax,\,\xmin,\,\p)$ is a
\Bs\ with associated spectral triple $(X,\,\cP_{2},\cR_{2})$. It is
easy to check that $\cR_{2}=\cR_{1}$. By Theorem \ref{t:pds}, $(X,
\xmax, \xmin, \p)$ is conjugate to a
Ver\v{s}ik transformation of an ordered Bratteli diagram. The
(unordered) Bratteli diagram defines the AF- algebra C$^*(\cR_2) =
$C$^*(\cR_1).$ Another way to see that $\cR_1$ is an AF groupoid is
to apply Theorem \ref{t:AFnest}.
Therefore, $\fA(\cP_{1}) $ is a
$\bZ$-analytic subalgebra of $C^{*}(\cR_{1})$ but $\fA(\cP_{1}) $ is
not standard $\bZ$-analytic .
\end{example}

\bibliographystyle{amsplain} \bibliography{opalg}

\end{document}